\documentclass[10pt,a4paper,oneside,twocolumn]{article}

\usepackage{geometry}
\usepackage{amsmath,amssymb}
\usepackage{graphicx}
\usepackage{psfrag}
\usepackage{color}
\usepackage[usenames,dvipsnames]{xcolor}
\usepackage{rotating}
\usepackage{nicefrac}
\usepackage{tikz}
\usetikzlibrary{shapes.callouts}
\usepackage{pgfplots}
\pgfplotsset{compat=1.8} 
\usepackage{listings}
\usepackage{tabularx}
\usepackage{mathabx}
\usepackage{stmaryrd}
\usepackage{units}
\usepackage[T1]{fontenc}
\usepackage[utf8]{inputenc}
\usepackage{changes}
\colorlet{Changes@Color}{red}

\usetikzlibrary{shapes,
                arrows,
                fit,
                backgrounds}
\definecolor{IPC}{RGB}{0,255,255}  
\definecolor{SPS}{RGB}{255,255,31} 
\definecolor{HIL}{RGB}{0,217,0}    

\newcommand{\norm}[1]{\left\lVert#1\right\rVert}
\newcommand{\R}{\mathbb{R}}
\newcommand{\X}{\mathcal{X}}
\newcommand{\Pe}{\mathcal{P}}
\newcommand{\N}{\mathbb{N}}
\newcommand{\A}{\mathcal{A}} 
\newcommand{\I}{\mathcal{I}} 
\newcommand{\U}{\mathcal{U}} 
\newcommand{\M}{\mathcal{M}} 
\newcommand{\T}{\mathcal{T}}

\newcommand\Toprule{\rule{0pt}{12px}}

\newcommand{\stkout}[1]{\ifmmode\text{\sout{\ensuremath{#1}}}\else\sout{#1}\fi}

\makeatletter
\newcommand{\xdashrightarrow}[2][]{\ext@arrow 0359\rightarrowfill@@{#1}{#2}}
\newcommand{\xdashleftarrow}[2][]{\ext@arrow 3095\leftarrowfill@@{#1}{#2}}
\def\rightarrowfill@@{\arrowfill@@\relax\relbar\rightarrow}
\def\leftarrowfill@@{\arrowfill@@\leftarrow\relbar\relax}
\def\arrowfill@@#1#2#3#4{%
  $\m@th\thickmuskip0mu\medmuskip\thickmuskip\thinmuskip\thickmuskip
   \relax#4#1
   \xleaders\hbox{$#4#2$}\hfill
   #3$%
}
\makeatother

\newtheorem{theorem}{Theorem}
\newtheorem{lem}[theorem]{Lemma}

\newtheorem{prop}[theorem]{Proposition}

\newcounter{myapproach}
\newcommand{\approach}[1]{\refstepcounter{myapproach} A\arabic{myapproach} }

\usepackage{titlesec}
\titleformat*{\section}{\large\bfseries}
\titleformat*{\subsection}{\normalsize\bfseries}
\titleformat*{\subsubsection}{\normalsize\bfseries}
\titleformat*{\paragraph}{\normalsize\bfseries}
\titleformat*{\subparagraph}{\normalsize\bfseries}

\titlespacing*{\section}{0pt}{1.5ex plus 1ex minus .2ex}{1.3ex plus .2ex}
\titlespacing*{\subsection}{0pt}{1.5ex plus 1ex minus .2ex}{1.3ex plus .2ex}

\begin{document}

\title{A Complexity Analysis of Event-Triggered Model Predictive Control on Industrial Hardware}

\author{{Patrik Simon Berner and Martin M\"onnigmann}%
\\Automatic Control and Systems Theory, Department of Mechanical Engineering, \\ Ruhr-Universit\"at Bochum, Germany,\\ \texttt{patrik-simon.berner@rub.de, martin.moennigmann@rub.de}%
}

\maketitle

\begin{abstract}
We implement a recently proposed event-triggered networked MPC approach on industrial hardware to analyze its practical relevance. 
There exist several alternatives for such an implementation that differ with respect to the distribution of computational load between local and central nodes, and with respect to network bandwidth requirements. 
These alternatives have been analyzed theoretically before, but when implemented it becomes evident that their usefulness cannot be predicted based on theoretical considerations alone. 
It is the purpose of the present paper to account for both practical and theoretical aspects in determining which alternative is most appropriate for an implementation on industrial hardware. 
The smallest possible bandwidth is known to result for a variant
in which only the active set of constraints is transmitted from the central to the local nodes. Since local nodes must determine the control law from the active set in this case, which requires matrix inversions, an unattractive computational cost results at first sight.   
Somewhat surprisingly, the computational cost scales practically linearly in the problem size when implemented.
We confirm this result with a more detailed theoretical complexity analysis than given in previous papers.   
All results are illustrated with data obtained with an implementation on industrial hardware components. 
\end{abstract}

\section{Introduction}
Model predictive control (MPC) is an established method for linear, multivariate systems with constraints.
MPC is computationally expensive, however.
MPC controllers are often too complex for an implementation on standard industrial hardware such as programmable logic controllers (PLC).

Various approaches to reducing the computational effort of MPC have been investigated. In \cite{Bemporad2002} and \cite{Seron2003}, parametric programming was used to solve the optimization problem offline.
Code generation tools (\cite{Zometa2013, Houska2011a, Mattingley2010}) produce a problem-tailored implementation of the MPC problem.
In \cite{Pannocchia2007} and \cite{Jost2015}, methods for an online simplification of the optimal control problem have been presented.
Event-based approaches, other than the one discussed here, have been proposed in~\cite{Eqtami2011, Varutti2009, JostSD2015, Incremona2017, Bernardini2012} and the references therein (see also the introduction in~\cite{Heemels2012}).

Tailored MPC implementations for hardware with limited resources such as PLCs have also been studied. A suboptimal parametric programming approach has been used in \cite{Palomo2011}. The degree of suboptimality was chosen so that the resulting controller satisfied resource restrictions of the PLC.
Parametric programming was also used in \cite{Rauova2011}, where a binary search tree is instrumental to an implementation with limited memory.
In \cite{Huyck2012}, the use of various quadratic program (QP) solvers for online MPC on PLC was investigated. The results were illustrated with a temperature control problem.
In \cite{Kufoalor2014} primal-dual first-order methods were used to speed up solving the optimization problem and the proposed method was tested with a hardware-in-the-loop (HIL) simulation.

We use an event-triggered MPC scheme~\cite{BernerP2016} that we claim to be particularly suitable for networked MPC.  
The central idea is as follows:
Solving a linear-quadratic optimal control problem (OCP) at the current state is usually assumed to provide the optimal control action for this particular state.
It is well-known, however, that this solution additionally provides an affine state feedback law that is optimal on a state space polytope~\cite{Bemporad2002,Seron2003}. 
As long as the system state remains in this polytope, the feedback law can be evaluated at very low computational effort.
Consequently, lean \textit{local nodes} such as PLCs, directly attached to controlled systems, suffice to evaluate affine feedback laws as long as the system stays in the polytope of validity. Upon leaving the polytope, a powerful \textit{central node} is called to solve an OCP and to provide a new feedback law and polytope.  
This setup obviously requires affine feedback laws (and linear inequalities that define polytopes) to be sent across the network.
Since bandwidth and computational resources are limited, it is crucial to investigate alternatives to sending the feedback law and polytope (i.e., real matrices and vectors) across the network. 
The present paper extends earlier theoretical results~\cite{BernerP2016} by new results on the computational and bandwidth requirements of the alternatives (Propositions~\ref{prop:compareTransmittedData}-\ref{prop:upperBoundKappa}), by details on the implementation (Section~\ref{sec:implementation}), and by a validation with an actual implementation on industrial hardware (Section~\ref{sec:results}). 

We introduce the event-triggered MPC scheme in Section~\ref{sec:triggerEventIntro}.
Section~\ref{sec:theoreticalPart} discusses the computational effort and bandwidth requirements of the event-triggered controller.
Section~\ref{sec:implementation} addresses the implementation on industrial hardware. 
Some results obtained with this implementation, which corroborate our theoretical complexity analysis, are discussed in Section~\ref{sec:results}.
A conclusion is drawn in Section~\ref{sec:conclusion}.

\subsubsection*{Notation}
Let $M_\A$ refer to the submatrix of the matrix $M$ with the row indices listed in the index set $\A$, where all index sets are assumed to be ordered.
A polytope $\Pe$ is defined to be the intersection of a finite number of halfspaces, i.e., for any polytope $\Pe$ there exist  $T \in \R^{q \times n}, \ d \in \R^{q}$ such that
\begin{equation}
        \Pe= \lbrace x \in \R^n \mid Tx \leq d \rbrace.
\end{equation}

\section{Event-triggered networked MPC}\label{sec:triggerEventIntro}

\subsection{Problem class}
We consider linear discrete-time MPC problems, where the optimal control problem
\begin{subequations}\label{eq:OCPformulation}
\begin{align}
  \scalebox{0.87}{$\min\limits_{\substack{u(k),\, k=0, \dots, N-1 \\ x(k),\, k= 1, \dots, N}} x(N)'$}&\scalebox{0.87}{$Px(N) + \sum\limits_{k=0}^{N-1} x(k)'Qx(k) +u(k)'Ru(k)$}\\
  \scalebox{0.87}{$\text{s.\,t.~~}x(k+1)~$}&\scalebox{0.87}{$= Ax(k)+Bu(k),\, k=0, \dots, N-1$} \label{eq:processModel}\\
  \scalebox{0.87}{$x(0)~$}&\scalebox{0.87}{$\text{given}$},\nonumber \\
  \scalebox{0.87}{$u(k)~$}&\scalebox{0.87}{$\in \U,\, k= 0, \dots, N-1$} \label{eq:Uconstraint}\\
  \scalebox{0.87}{$x(k)~$}&\scalebox{0.87}{$\in \X,\, k= 0, \dots, N-1$} \label{eq:Xconstraint}\\
  \scalebox{0.87}{$x(N)~$}&\scalebox{0.87}{$\in \T$}, \label{eq:terminalConstraint}
\end{align}%
\end{subequations}%
with state variables $x \in \R^n$ and input variables $u\in\R^m$, horizon $N \in \N >1$, weighting matrices $P\in\R^{n\times n}$, $P\succ 0$, $Q\in\R^{n\times n}$, $Q\succeq 0$ and $R\in\R^{m\times m}$, $R\succ 0$ is solved.
The system $(A,B)$ with matrices $A\in\R^{n\times n}$ and $B\in\R^{n\times m}$ is assumed to be stabilizable.
States and inputs are subject to constraints \eqref{eq:Uconstraint}, \eqref{eq:Xconstraint}, \eqref{eq:terminalConstraint} with compact polytopes $\U\subset\R^m$, $\X\subset\R^n$ and $\T\subseteq\X$ that contain the origin in their interiors.

The MPC problem \eqref{eq:OCPformulation} can be stated as a quadratic program (QP) of the form~\cite[Chap. 3]{Maciejowski2002}:
\begin{equation}\label{eq:condensedMPC}
	\min_U \frac{1}{2}U'HU + x'FU \quad \text{s.t. } GU-Ex \leq w
\end{equation}
with $U= (u(0)^\prime, \dots, u(N-1)^\prime)^\prime \in \R^{mN}$, $F \in \R^{n \times mN}$, $H \in \R^{mN \times mN}$, $G \in \R^{q \times mN}$, $w \in \R^{q}$, $E \in \R^{q \times n}$, where $\mathcal{Q} =  \left\lbrace 1, \hdots, q \right\rbrace $ and $q$ is the number of constraints in~\eqref{eq:OCPformulation} and~\eqref{eq:condensedMPC}.
$H$ is positive definite since $Q\succeq 0 $ and $R\succ 0$~\cite[p. 76]{Maciejowski2002}.
The set of all feasible states is denoted by $\X_f$ and assumed to be nonempty. Since $H\succ 0$, there exists a unique optimal solution for every $x\in\X_f$ that we denote $U^{\star}(x)$.
There exist polytopes $\Pe_1, ..., \Pe_p$ and affine control laws \mbox{$K_ix+b_i$}, $K_i \in\R^{mN \times n}$, $b_i \in\R^{mN}$, $i=1,...,p$, such that  $\cup_i^p \Pe_i = \X_f$ and
\begin{equation}\label{eq:PWAsolution}
	U^{\star}(x)= \left\{\begin{array}{cl} K_1x+b_1 \text{, if } x\in\Pe_1,\\ \vdots \\K_px+b_p \text{, if } x\in\Pe_p, \end{array}\right.
\end{equation}
and $U^{\star}(x)$ is a continuous function on $\X_f$~\cite{Bemporad2002,Seron2003}.
We stress that we exploit the structure of the explicit solution~\eqref{eq:PWAsolution} without ever calculating it entirely.

\subsection{Event-triggered controller}
Equation~\eqref{eq:PWAsolution} suggests the simple networked MPC setup that was already briefly described in the introduction:
  Solving~\eqref{eq:OCPformulation} or~\eqref{eq:condensedMPC} for the current state $x(k)\in\R^n$ is usually understood to result in the optimal input signal for this particular point $x(k)$ in the state space. 
  It follows from~\eqref{eq:PWAsolution}, however, that the solution to~\eqref{eq:condensedMPC} for a point $x\in\X_f$ not only provides the optimal input signal at this point, but an affine law $K_ix+b_i$ and its polytope of validity $\mathcal{P}_i$ (see Lemma~\ref{lem:OptControllerGainAndPolytope} below).  
 For all $x\in\mathcal{P}_i$, this affine law yields the same optimal input as solving the optimal control problem~\eqref{eq:OCPformulation} or~\eqref{eq:condensedMPC}. 
Since the computational effort for evaluating $K_i x+ b_i$ and checking whether $x\in\mathcal{P}_i$ is much smaller than that for solving~\eqref{eq:OCPformulation} or \eqref{eq:condensedMPC}, the optimal solution can be determined on a lean local node such as a PLC 
with $K_i x+ b_i$ and $\mathcal{P}_i$.
A computationally powerful central node (the industrial PC in the particular setup investigated here) is only required to solve the OCP whenever the current polytope is left.

The calculations necessary to obtain the affine law $K_i x+ b_i$ and its polytope of validity $\mathcal{P}_i$ from $U^\star(x)$ are summarized in the following Lemma for ease of reference. 
\begin{lem} 
\label{lem:OptControllerGainAndPolytope} 
Let $U^\star(x)$ be the solution to \eqref{eq:condensedMPC} for an arbitrary $x\in\mathcal{X}_f$, and let $\A$ and $\I$ refer to the sets of active and inactive constraints
\begin{equation}
\label{eq:IndexSetsAandI} 
 \begin{split}
   \mathcal{A} &= \left\{i\in \left\{1,...,q\right\} \left|G_{i} \,U^\star(x)\!-E_{i} x= w_{i} \right.\right\},
   \\
   \mathcal{I} &= \left\{i\in \left\{1,...,q\right\}\left|G_{i} \,U^\star(x)\!-E_{i} x< w_{i} \right.\right\}.
 \end{split}
\end{equation}
Let $\M=\{1, \dots, m\}$ and assume the matrix $G_\A$ has full row rank. Then%
\begin{equation}\label{eq:MatricesForAffineLawAndPolytope}
    \begin{aligned}
    \scalebox{0.93}{$K_\M^{\star}$} &\scalebox{0.93}{$=H_\M^{-1}(G_{\A})'(G_{\A}H^{-1}(G_{\A})')^{-1}S_{\A}-H_\M^{-1}F'$},\\
    \scalebox{0.93}{$b_\M^{\star}$} &\scalebox{0.93}{$=H_\M^{-1}(G_{\A})'(G_{\A}H^{-1}(G_{\A})')^{-1}w_{\A}$},\\
    \scalebox{0.93}{$T^{\star}$} &\scalebox{0.93}{$= \begin{pmatrix} G_{\I}H^{-1}(G_{\A})'(G_{\A}H^{-1}(G_{\A})')^{-1}S_{\A}-S_{\I} \\    (G_{\A}H^{-1}(G_{\A})')^{-1}S_{\A}    \end{pmatrix}$},\\
    \scalebox{0.93}{$d^{\star}$} &\scalebox{0.93}{$= -\begin{pmatrix} G_{\I}H^{-1}(G_{\A})'(G_{\A}H^{-1}(G_{\A})')^{-1}w_{\A}-w_{\I}\\ (G_{\A}H^{-1}(G_{\A})')^{-1}w_{\A}\end{pmatrix},$}
\end{aligned}
\end{equation} 
with $S = E+G H^{-1}F^\prime$, $S\in\R^{q\times n}$
define the optimal control law and the polytope 
\begin{subequations}\label{eq:OptControllerGainAndPolytope}
\begin{align}
 u^\star(x) &= K_\M^\star x + b_\M^\star \text{ } \forall \text{ } x\in\Pe^\star,\label{eq:OptControllerGain}\\
  \Pe^\star &= \{ x \in \R^n \,|\, T^\star x \leq d^\star \}. \label{eq:OptControllerPolytope}
\end{align}
\end{subequations}
\end{lem}
A concise proof of Lemma~\ref{lem:OptControllerGainAndPolytope}, which is an immediate corollary to the results in~\cite{Bemporad2002}, is given in~\cite{JostSD2015}. 

For later use we introduce $q_\A= |\A|$ and 
\begin{equation}\label{eq:triggeringCondition}
		e(x)= \left\{\begin{array}{ll}0 &\text{ if } x\in\Pe^\star, \\ 1 &\text{ otherwise.} \end{array}\right. 
\end{equation}
The function $e(x)$ is used to check whether the control law~\eqref{eq:OptControllerGain} can be reused to compute the control input for the current state $x(k)$ ($e(x(k))= 0$), or a new optimal control law has to be calculated ($e(x(k))= 1$). 
We introduce the term \textit{trigger event}, or in short \textit{event} to describe any instant for which $e(x(k))=1$.
The inequalities in~\eqref{eq:OptControllerPolytope} are validated row by row in the implementation. If all inequalities hold, no event occurred (i.e., $e(k)=0$).

\subsection{Event-triggered networked MPC}
\begin{figure}[tb]
\vspace*{0px}
\centering
	\psfrag{r}{\hspace{-38px}\parbox[t]{67px}{\vspace{-32px} \scriptsize 
	$\text{\textbf{if}}\, e\!=\!0$\newline
	\phantom{+}$u^\star \!\!:=\! K_\M^\star x \!+\! b_\M^\star$\newline
	\textbf{elseif $e\!=\!1$}\newline
	\phantom{+}\textit{request feedback}\newline
	\phantom{+} \textit{law and polytope}\newline
	\phantom{+}$u^\star \!\!:=\! K_\M^\star x \!+\! b_\M^\star$\newline 
	\textbf{end}
	}}
	\psfrag{s}{\hspace{-12px}\parbox[t]{30px}{\vspace{-15px}controlled \phantom{              }system} }
	\psfrag{u}{\!\small$u^\star$}
	\psfrag{x}{\small$x$}
	\centering
	\includegraphics[width=0.95\linewidth]{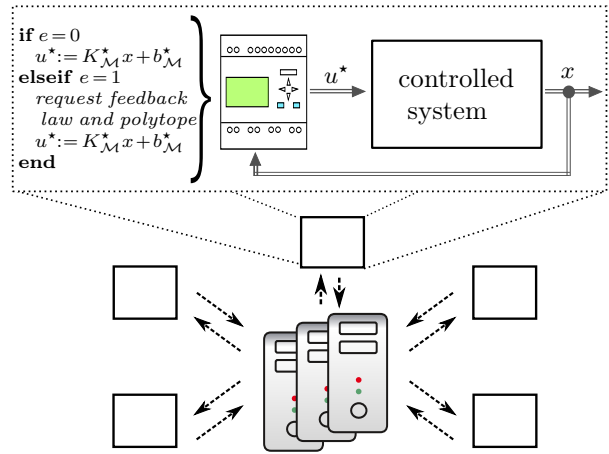}
	\caption{Structure of the event-triggered MPC. Each box represents a local node that controls an attached system in closed-loop control (see dotted box for an example). Dashed lines indicate network connections used on demand. The feedback loop in the local setup (double lines) is permanent. The local node evaluates the locally stored feedback law~\eqref{eq:OptControllerGain} until a trigger event occurs (i.e., until $e(x)=1$ in~\eqref{eq:triggeringCondition}). If $e(x)=1$, the local feedback law is replaced by a new law requested from the central node.}
    \label{fig:setup}
    \vspace*{0px}
\end{figure}

Consider the simple star-shaped network architecture depicted in Fig.~\ref{fig:setup},  where multiple  local control nodes are connected to a single central node via network connections.
Throughout this paper we call this structure an \textit{event-triggered networked model predictive controller} in contrast to \textit{classical MPC}, where the OCP is solved and the control input is applied to the system in every time step.
\\
In case of an event, a control law and polytope must be sent across the network. 
The simplest way of doing so is to send matrices $K_\M^\star, b_\M^\star, T^\star$ and $d^\star$, that represent the control law~\eqref{eq:OptControllerGain} and polytope~\eqref{eq:OptControllerPolytope}.
However, equations~(\ref{eq:IndexSetsAandI}-\ref{eq:OptControllerPolytope}) suggest to distribute computations between local and central node such that the amount of transmitted data may be reduced.
It is our aim to identify the approach that results in the best trade-off between transmitted data and computations on the local node, and to propose an implementation scheme for this particular variant.

%

\section{Computational- and network load in event-triggered networked MPC}
\label{sec:theoreticalPart}

\subsection{Measuring transmitted data and computational cost}
The remainder of the paper is based on the following assumptions and hypotheses  H1-H3:
(H1) Bandwidth requirements can be determined by counting the number of bits transferred per event, where real numbers are transmitted using IEEE 754 half precision~\cite{IEEE754} ($\lambda= 16\frac{\text{bit}}{\text{real}}$ below). 
(H2) Computational effort can be measured by counting floating point operations, where multiplications and summations cost one floating point operation (flop), and divisions cost ten flops~\cite{Ueberhuber1997};
the cost of some recurring operations is stated in Table~\ref{tab:basicops} \cite[chap. 1.1.15]{Golub2013} \cite[p.12]{Farebrother1988}.
(H3) The central node is assumed to  carry out all computations instantaneously.

Let the term \textit{computational load on the local node} refer to the number of flops necessary to determine the control law and polytope~\eqref{eq:OptControllerGainAndPolytope}. 
Since all variants discussed below have the evaluation of the control law~\eqref{eq:OptControllerGain} and condition~\eqref{eq:triggeringCondition} in common, we neglect the flops required for~\eqref{eq:OptControllerGain} and~\eqref{eq:triggeringCondition}.

\begin{table}[h]
\renewcommand{\arraystretch}{1.3}
\caption{Computational effort for some matrix operations}
\vspace*{0px}
  \centering
  \begin{normalsize}
  \begin{tabular}{l|l|l|l|l}
         op.   & $M + \tilde{M}$ & $c M$ & $M V$ 	 & $M^{-1}$\\ \hline 
\Toprule flops & $mn$			& $m n$ & $ml(2n-1)$ & $\frac{2n^3+18n^2+10n}{3}$
  \end{tabular}\\
  \end{normalsize}
\Toprule  {\begin{scriptsize}\textit{For any $M,\tilde{M} \in \R^{m \times n}$, $V \in \R^{n \times l}$, $c \in \R$ and $m, n, l\in \N$. Existence of the inverse and $m=n$ are assumed in the last column.}\end{scriptsize}}
  \label{tab:basicops}
\end{table}%

\subsection{Trade-off between transmitted data and local calculations}\label{subsec:Tradeoff}
A control law and polytope must be sent across the network whenever an event~\eqref{eq:triggeringCondition} occurs.
It is evident from~\eqref{eq:OptControllerGainAndPolytope} that the matrices $K_\M^\star, b_\M^\star, T^\star$ and $d^\star$ define the control law and polytope.
It suffices to send much less data, however. This is stated more precisely in Lemma~\ref{lem:OptionsForTransmittedData}, which was proved in~\cite{BernerP2016}.
\begin{lem}
\label{lem:OptionsForTransmittedData}
Let $x \in \X_f$ be arbitrary and assume the optimal control problem~\eqref{eq:condensedMPC} has been solved for this initial condition. Then the feedback law~\eqref{eq:OptControllerGain} and its polytope of validity~\eqref{eq:OptControllerPolytope} are determined by any of the following:
\begin{subequations}
  \begin{align}
  & \A(x),
  \label{eq:TransferActiveSet}
  \\
  & U^{\star}(x),
  \label{eq:TransferInputSequence}
  \\
  & K_{\M}^{\star},\, b_{\M}^{\star},\, T^{\star},\, d^{\star}.
  \label{eq:TransferMatricesAndVectorsKbTd}
  \end{align}
\end{subequations}
\end{lem}

Equation~\eqref{eq:TransferInputSequence} may tempt the reader to think we use the input signals $u^\star(0), \dots, u^\star(N-1)$ of $U^\star$, 
which is open-loop but in general not closed-loop optimal. 
We stress again this is not the case, but we use the feedback law \eqref{eq:OptControllerGain} and polytope~\eqref{eq:OptControllerPolytope}, 
because this feedback law is closed-loop optimal on the polytope. 
Because we determine a new law and polytope whenever the current polytope is left, 
the same closed-loop optimal sequence as solving~\eqref{eq:OCPformulation} in every time step results.
Therefore, the same closed-loop optimal system behavior is obtained as solving~\eqref{eq:OCPformulation} or~\eqref{eq:condensedMPC} in every time step~\cite{JostSD2015}. 
The claim about~\eqref{eq:TransferInputSequence} in Lemma~\ref{lem:OptionsForTransmittedData} holds, because $U^\star$ defines $\mathcal{A}(x)$ according to~\eqref{eq:IndexSetsAandI} and $\mathcal{A}$ defines $K^\star_{\mathcal{M}}$, $b^\star_{\mathcal{M}}$, $T\star$, $d^\star$ and thus the law and polytope according to \eqref{eq:MatricesForAffineLawAndPolytope}.

Lemma~\ref{lem:OptionsForTransmittedData} implies the following alternative approaches A1-A4 for transmitting the optimal control law~\eqref{eq:OptControllerGain} and its polytope~\eqref{eq:OptControllerPolytope}.

\begin{itemize}
	\item[\approach{A1}\label{appr:activeset}] Send $q$ bits 
		\begin{equation}\label{eq:bitwiseRepresentationGamma}
		  \gamma_i= \left\{\begin{array}{cl} 0 & \mbox{if }G_iU^{\star}+E_ix < w_i\\ 1 & \mbox{if } G_iU^{\star}+E_ix = w_i \end{array}\right. ,
		\end{equation} 
		$i= 1, \dots, q$, from the central to the local node, which is equivalent to sending $\mathcal{A}$ and $\mathcal{I}$ according to~\eqref{eq:IndexSetsAandI}. 
		Evaluate~\eqref{eq:MatricesForAffineLawAndPolytope} on the local node to determine $K_{\M}^{\star}$, $b_{\M}^{\star}$, $T^\star$ and $d^\star$.
	\item[\approach{A2}\label{appr:sendInverse}] Determine 
		\begin{equation}
		\label{eq:Phi}
		\Phi = (G_{\A}H^{-1}(G_{\A})')^{-1}
		\end{equation}
		on the central node and send it in addition to and redundantly to~\eqref{eq:bitwiseRepresentationGamma}. 
		Evaluate ~\eqref{eq:MatricesForAffineLawAndPolytope} to determine $K_{\M}^{\star}$, $b_{\M}^{\star}$, $T^\star$ and $d^\star$ as in A1 on the local node.
		Use the inverse $\Phi$ instead of determining it on the local node.
	 \item[\approach{A3}\label{appr:optsolution}] Send $U^\star$ from the central to the local node.  
	 Evaluate~\eqref{eq:IndexSetsAandI} to determine $\mathcal{A}$ and $\mathcal{I}$, then evaluate~\eqref{eq:MatricesForAffineLawAndPolytope}  to determine 
	 $K_{\M}^{\star}$, $b_{\M}^{\star}$, $T^\star$ and $d^\star$ on the local node.
 \item[\approach{A4}\label{appr:firstM}] Determine $K_{\M}^{\star}$, $b_{\M}^{\star}$,  $T^\star$, $d^\star$ on the central node and send them to the local node.
\end{itemize}
These four approaches obviously differ with respect to the amount of data sent across the network and the computational effort on the local node. 
More precisely, we have the following results, which are taken from~\cite{BernerP2016}.
\begin{lem}\label{lem:approachesAndEffort}
Assume H1-H4 hold. Furthermore,  
assume~\eqref{eq:OCPformulation} or equivalently~\eqref{eq:condensedMPC} has been solved for an $x\in\mathcal{X}_f$
on the central node.
Then the number of bits that need to be transmitted in approaches A1-A4 and the number of flops required on the local node for approaches A1-A4 are as follows: 

\vspace{1mm}
\footnotesize
\begin{tabular}[tb]{lll}
	& bits &  flops \\
	 \hline
	 \\
	A1 & $q$ &
	 $(mNq_{\A}+mn+qq_{\A})(2mN-1)$\\
	&& $+(qn+q+m+mn)(2q_{\A}-1)$\\
	&& $+\frac{2}{3}q_{\A}^3+6q_{\A}^2+\frac{7}{3}q_{\A}$\\
	&& $+qn+q-q_{\A}n+mn$
	   \\
	   \\
	A2 & 	$\frac{\lambda \left(q_{\A}^2+q_{\A}\right)}{2}+q$ &
	$(mNq_{\A}\!+\!mn\!+\!qq_{\A}\!-\!q_{\A}^2)(2mN\!\!-\!\!1)$\\
	&& $+(qn+q+m+mn)(2q_{\A}-1)$\\
	&& $+qn+q-q_{\A}n-q_{\A}+mn$
	   \\
	   \\
	A3 & $\lambda mN$ &
	$(mNq_{\A}\!+\!mn\!+\!qq_{\A}\!+\!q)(2mN\!\!-\!\!1)$\\
	&& $+(qn+q+m+mn)(2q_{\A}-1)$\\
	&& $+\frac{2}{3}q_{\A}^3+6q_{\A}^2+\frac{7}{3}q_{\A}+3qn$\\
	&& $+2q-q_{\A}n+mn$
	   \\
	   \\
	A4 & $\lambda \left(mn + m \right.$ & $0$\\	
	& $\left.+ qn +q\right)$  &
\end{tabular}
\normalsize
\vspace{4px}
\end{lem}

Based on the technical results stated in Lemma~\ref{lem:approachesAndEffort}, it is now easy to compare approaches A1-A4.
We compare the amount of transmitted data first. 
\begin{prop}\label{prop:compareTransmittedData}
Let all assumptions be as in Lemma~\ref{lem:OptionsForTransmittedData}. 
The following statements hold for the number of transmitted bits in approaches~\ref{appr:activeset}-\ref{appr:firstM}:
\begin{itemize}
\item[(i)] Approach~\ref{appr:activeset} requires as much as, or less data than approach~\ref{appr:sendInverse} to be sent across the network.
\end{itemize}
 If we additionally assume all constraints (\ref{eq:Uconstraint},~\ref{eq:Xconstraint},~\ref{eq:terminalConstraint}) to be box constraints 
 (i.e., if $q=2mN+2nN+2n$),
  the following statements hold:
 \begin{itemize}
\item[(ii)] Approaches~\ref{appr:activeset}-\ref{appr:optsolution} never require more data to be transmitted than approach~\ref{appr:firstM}.
\item[(iii)] Approach~\ref{appr:optsolution} requires more data to be sent across the network than approach~\ref{appr:activeset}, if
\begin{align}\label{eq:compApprABbits}
	\frac{\lambda-2}{3} > \frac{n}{m}
\end{align}
where $\lambda$ is the number of bits per real number.
\item[(iv)] Approach~\ref{appr:optsolution} requires more data to be sent across the network than approach~\ref{appr:sendInverse}, if 
\begin{align}\label{eq:compAppr32bits}
\frac{\lambda-2}{3} > \frac{n}{m} + \frac{\lambda (q_{\A}(q_{\A}+1))}{6mN}.
\end{align}
\end{itemize}
\end{prop}

\textit{Proof:}
Statement~(i) is obvious, since \ref{appr:activeset} requires to transmit $\gamma$ and approach~\ref{appr:sendInverse} requires to transmit $\gamma$ and $\Phi$.
Statement (ii) can be shown by substituting $q=2mN+2nN+2n$ (box constraints) into the number of transmitted bits stated in Lemma~\ref{lem:approachesAndEffort} and subtracting the number of bits for~\ref{appr:activeset}-\ref{appr:optsolution} from those for~\ref{appr:firstM}. Since the resulting expressions are always positive, statement (ii) holds.
In order to prove (iii) we need to show that~\eqref{eq:compApprABbits} implies\
  \begin{align}\label{eq:statement3Lemma3}
    \lambda mN > q .
  \end{align}
For $q=2mN+2nN+2n$ (box constraints), this is equivalent to
\begin{align}\label{eq:statement3Lemma3box}
	 \lambda > 2+2\frac{n}{m}(1+\frac{1}{N}).
\end{align}
It is easy to show $N>1$ implies $2\frac{n}{m}\left( 1+ \frac{1}{N} \right) \leq 3 \frac{n}{m}$. Therefore,
\begin{align}\label{eq:statement3Lemma3boxN}
	 \lambda > 2+3\frac{n}{m}
\end{align}
implies~\eqref{eq:statement3Lemma3box}, i.e., \eqref{eq:statement3Lemma3} for box constraints. 
Rearranging~\eqref{eq:statement3Lemma3boxN} yields the desired condition~\eqref{eq:compApprABbits}.
Part~(iv) can be shown in the same way as~(iii).
$\hfill \blacksquare$
Since states and inputs are constrained by an upper and lower bound in many technical systems, the introduction of box constraints in part (ii) - (iv) of Proposition~\ref{prop:compareTransmittedData} is not a severe restriction.

Proposition~\ref{prop:compareTransmittedData} reveals the following partial order of the approaches with respect to the number of transmitted bits: Approach~\ref{appr:activeset} requires the smallest number of transmitted bits, followed by approaches~\ref{appr:sendInverse} and~\ref{appr:firstM}. Approach~\ref{appr:optsolution} results in less transmitted data than~\ref{appr:activeset} if inequality~\eqref{eq:compApprABbits} is \textit{not} fulfilled and less transmitted data than approach~\ref{appr:sendInverse} if inequality~\eqref{eq:compAppr32bits} is \textit{not} fulfilled. 
It remains to compare the computational effort on the local node based on Lemma~\ref{lem:OptionsForTransmittedData}.
\begin{prop}\label{prop:compareComputations}
Let all assumptions be as in Lemma~\ref{lem:OptionsForTransmittedData}. The following statements hold for the computational effort on the local node. 
\begin{itemize}
\item[(i)] The number of flops in approach~\ref{appr:sendInverse} never exceeds those needed in approach~\ref{appr:activeset}.
\item[(ii)] Approach~\ref{appr:optsolution} always requires more computations than approach~\ref{appr:activeset}.
\item[(iii)] Approach~\ref{appr:firstM} requires no  operations on the local node to determine the matrices for the control law and polytope.
\end{itemize}
\end{prop}
\textit{Proof:}
Statement (iii) is obvious from the zero in the last row and column of the table in Lemma~\ref{lem:approachesAndEffort}. 
Statements (i) and (ii) can be shown by subtracting the number of flops given in the third column of the table in Lemma~\ref{lem:approachesAndEffort} for the respective pair of approaches. Specifically, this difference reads
$\frac{2}{3}q_{\A}^3+q_{\A}^2(2mN+5)+\frac{10}{3}q_{\A}$ for approach~\ref{appr:activeset}. 
Since this expression is positive, approach~\ref{appr:activeset} never requires less operations than approach~\ref{appr:sendInverse}. Statement (ii) can be shown in the same way.~$\hfill \blacksquare$

None of the approaches~\ref{appr:activeset}-\ref{appr:firstM} results in both, the smallest number of transmitted bits and the smallest computational effort. 
Since the number of transmitted bits in approach~\ref{appr:sendInverse} depends on $q_\A$ (see~\eqref{eq:compAppr32bits}), and since approach~\ref{appr:optsolution} depends on~\eqref{eq:compApprABbits}, the best trade-off results for transmitting the active set of constraints in approach~\ref{appr:activeset}.
We examine \ref{appr:activeset} in more detail in the remainder of the paper.

\subsection{Computational effort of approach~\ref{appr:activeset}}
While only $q_\A$ bits need to be transmitted in approach~\ref{appr:activeset}, additional computations must be carried out on the local node. Most importantly, the inverse of $G_{\A}H^{-1}(G_\A)'\in\R^{q_\A\times q_\A}$ has to be determined on the local node, which requires on the order of $q_\A^3$ operations (see, e.g., \cite{Farebrother1988}, p.12). However, the absolute number of operations for this inversion always remains small compared to the number of operations necessary for evaluating~\eqref{eq:MatricesForAffineLawAndPolytope}, even though the latter scales linearly only in $q_\A$. 
This statement is made more precise in the following proposition for the special case of box constraints. Note that the proposition is conservative in the sense that the computational effort for the inversion is even smaller in practice (see Section~\ref{sec:implementation}). 
\begin{prop}\label{prop:upperBoundKappa}
Let $\eta_{\rm inv}(q_\A) $ and $\eta_{\rm mat}(q_\A)$ refer to the number of flops necessary to carry out the matrix inversion~\eqref{eq:Phi} and and all other matrix operations in~\eqref{eq:MatricesForAffineLawAndPolytope}, respectively.
Assume the constraints~\eqref{eq:Uconstraint}-\eqref{eq:terminalConstraint} to be box constraints, which implies $q=2mN+2nN+2n$, and assume $N>1$. 
Then  
\begin{align}\label{eq:kappa}
	\frac{\eta_{\rm inv}(q_\A)}{\eta_{\rm mat}(q_\A)} \leq \frac{18}{79}.
\end{align}
\end{prop}

\textit{Proof:}
The inversion in~\eqref{eq:MatricesForAffineLawAndPolytope} requires
\begin{align*}
  \eta_{\rm inv}(q_\A) &= \frac{1}{3}(2q_{\A}^3 + 18 q_{\A}^2 + 10q_\A)
\end{align*}
flops according to~\cite[p.12]{Farebrother1988}. Subtracting this figure from the overall number of operations for~\ref{appr:activeset} as stated in Lemma~\ref{lem:approachesAndEffort} yields
\begin{align*}
  	\eta_{\rm mat}(q_\A) =& \alpha q_\A + \beta \text{, with}\\
  	\alpha =& (mN+q)(2mN-1)-n-1+2m+2mn\\ &+2qn+2q,\\
  	\beta =& -m+mn(2mN-1),
\end{align*}
which equals the number of flops for all matrix operations but the inversion in~\ref{appr:activeset}.
It is easy to show that $m>0$, $N>1$, $q>0$ and $n>0$ imply $\alpha>0$ and $\beta>0$.  
We now show that if
\begin{align}\label{eq:upperBoundOps}
  \frac{ \eta_{\rm inv}(q_\A)}{\eta_{\rm mat}(q_\A)} \leq \frac{18}{79}
\end{align}
holds for a $\bar{q}_{\A}$, it holds for all $q_\A< \bar{q}_\A$. 
This can be seen by extending $q_\A\in\N$ to a real variable $q_\A\in\R$ and showing $\eta_{\rm inv}/\eta_{\rm mat}$ is strictly increasing by proving $\left. \frac{d}{d q_\A}\! \left(\frac{\eta_{\rm inv}(q_\A)}{\eta_{\rm mat}(q_\A)}\right)\! > 0 \right.$. 
Applying the quotient rule yields a fraction with numerator
\begin{align}\label{eq:nominatorUpperBound}
\!(2q_{\A}^2 \! + \! 12q_{\A} \! + \! \frac{10}{3})(\alpha q_{\A} \! + \! \beta) \! - \! (\frac{2}{3} q_{\A}^2 \! + \! 6q_{\A} \! + \! \frac{10}{3}) \alpha q_{\A}
\end{align}
and denominator $\eta_{\rm mat}^2(q_\A)$. Since the denominator is positive, it suffices to show~\eqref{eq:nominatorUpperBound} is positive, which can equivalently be stated as 
\begin{align*}
  (\frac{4}{3} q_{\A}^2 + 6q_{\A}) \alpha q_{\A} + (2q_{\A}^2 + 12q_{\A} + \frac{10}{3}) \beta > 0.
\end{align*}
This condition obviously holds for $\alpha >0$, $\beta > 0$ and $q_\A > 0$.
It remains to be shown that~\eqref{eq:upperBoundOps} is positive for the largest possible number of concurrently active constraints $\bar{q}_\A$. 
This number reads $q_\A= \bar{q}_\A= mN$ for 
box constraints with $q=2mN+2nN+2n$. All combinations of $N>1$, $n>0$ and $m>0$ are covered by the three cases
\begin{align*}
\mbox{(i)}\, & N= 2,\, n= 1,\, m\ge 1, \\
\mbox{(ii)}\, & N= 2,\, n\ge 2,\, m\ge 1,\, \\
\mbox{(iii)}\, & N\ge 3,\, n\ge 1,\, m\ge 1, 
\end{align*}
where (i) and (ii) together cover $N= 2$, $n> 0$, $m> 0$. 
Case (i), i.e., substituting $N= 2$ and $n= 1$ into~\eqref{eq:upperBoundOps} and rearranging yields the condition
$1328 m^2-1368 m+40 \ge 0$, which obviously holds for all $m\ge 1$. 
Case (ii), i.e., substituting $N= 2$ into~\eqref{eq:upperBoundOps} and rearranging yields the condition
\begin{align}\label{eq:ProofLemma5Helper1}
 \!\!\!\! 1328m^2\!\! +\! 3888mn \! +\! 1296n^2 \!\! -\! 5256m\! +\! 486n\! -\! 1742\ge 0.
\end{align}
Now $n\ge 2$ implies 
$3888mn-5256m\ge 0$ and $1296n^2+486n-1742 \ge 0$ for all $m$. Adding these two inequalities to $1328m^2\ge 0$, which also holds for all $m$, shows that condition~\eqref{eq:ProofLemma5Helper1} holds. 
Case (iii) can be shown by noting that~\eqref{eq:upperBoundOps} can equivalently be stated as
\begin{align}\label{eq:ProofLemma5Helper2}
 & 166 N^3 m^2+216 N^3 m n+432 N^2 m n+216 N^2 n^2
 \\\nonumber
 +&108 N^2 n+216 N m n+216 N n^2+108 N m+54 N n
 \\
 \ge &1368 N^2 m+844 N+54 n+54
 \nonumber
\end{align}
for arbitrary $N$, $m$ and $n$. It is easy to see that $N\ge 3$ and $m\ge 1$ imply the first three terms on the left hand side of~\eqref{eq:ProofLemma5Helper2} are larger than the first term on the right hand side for all $\left. n\ge 1 \right.$. Similarly, the 4th and 5th term on the left hand side of~\eqref{eq:ProofLemma5Helper2} together dominate the second term on the right hand side, and the last term on the left hand side dominates the remaining two terms on the right hand side. Since the remaining terms on the left hand side are positive, condition~\eqref{eq:ProofLemma5Helper2} holds. 
$\hfill \blacksquare$

\section{Implementation on industrial hardware}\label{sec:implementation}
We use a SIMATIC IPC847C industrial computer (IPC for short) with an Intel Core i7-610E CPU and 8GB memory\footnote{\label{note:Siemens}See \texttt{http://www.siemens.com} for more information about the IPC, PLC and I/O modules (checked on October 7, 2016).} as the central node. 
The IPC runs a QNX Neutrino 6.5.0 real-time operating system and is programmed in ISO C++.
The local node consists of a Siemens S7-400 PLC equipped with 30MB memory, an SM 331 module, and an SM 332 analog I/O module\footnotemark[2]. We use statement lists (STL) to program the PLC. 
The system that is controlled by the PLC is simulated on a dSpace\footnote{
	See \texttt{https://www.dspace.com} for more information  
	(checked on October 7, 2016).
} HIL system. System states and inputs are exchanged via analog I/Os. 
We refer to this equipment as \textit{HIL setup} for brevity. 

Computation times are measured on the PLC by setting a digital output to high every time the computational task starts, and to low when it finishes. 
Closer inspection of these measurements reveals that the digital I/O modules of the PLC have a small random delay.
We repeat measurements ten times to reduce the influence of this delay. This results in standard deviations of less than 3\% for all reported data.

The communication between the nodes is performed by a PROFIBUS network, but we stress that any other bus system can be used instead. The network is only used when the local node triggers the central node to calculate a new affine control law and polytope.

\subsection{Implementation details}\label{subsec:implementationDetails}
We use a simple code generator to generate the required STL code from the problem \eqref{eq:OCPformulation} in Matlab. The code generator precomputes matrices $H^{-1},G,S,F',w$ offline.
The evaluation of~\eqref{eq:MatricesForAffineLawAndPolytope} exploits common subexpressions. 
The inversion of $G_\A H^{-1} (G_\A)'$ is replaced by a LU decomposition with pivoting and forward-backward substitutions. 
The code generator also generates a C++ header file for the IPC which provides matrices of the QP~\eqref{eq:condensedMPC}.  We use OOQP (see \cite{Gertz2003}) to solve QPs on the IPC.

Since data transmission and calculations~\eqref{eq:condensedMPC},~\eqref{eq:MatricesForAffineLawAndPolytope} may be necessary to obtain the control input $u$, a delay between measuring the state and applying the control input may occur. 
However, the delay only occurs if $e(k)=1$ in~\eqref{eq:triggeringCondition}.
Consequently, depending on whether $e(k)=0$ or $e(k)=1$ in~\eqref{eq:triggeringCondition}, the computation of the control input differs. We explain these two cases in more detail.
Let $x(k)$ and $x(k+1)$ refer to the current and subsequent state of the closed-loop system, and let $\Pe^\star_k$ be a polytope with $x(k) \in \Pe^\star_k$. While $x(k)$ is known from a measurement, $x(k+1)$ is predicted at sampling interval $k$ with the system model~\eqref{eq:processModel}. 
Either one of the following two cases applies in every time step:

(i) If $x(k+1) \in \Pe^\star_k$ , the control law and polytope in the subsequent time step $k+1$ remain the same as in time step $k$.
The control input $u(k+1)$ is calculated at the beginning of sampling interval $k+1$ on the local node and is applied to the system. In this case \textit{no} communication with the central node is necessary.

\begin{figure}[t]
	\centering
\centering
\begin{tikzpicture}[auto, font=\scriptsize,
                    cloud/.style={ellipse,
                                  draw=black,
                                  fill=white,
                                  minimum height=8mm,
                                  text centered,
                                  text width=21mm,
                                  inner sep=1pt,
                                 },
                    block/.style={rectangle,
                                  draw=black,
                                  fill=white,
                                  rounded corners=1mm,
                                  minimum height=11mm,
                                  text centered,
                                  text width=27mm,
                                 },
                    decis/.style={diamond,
                                  aspect=1.5,
                                  draw=black,
                                  thick,
                                  fill=white,
                                  text width=16mm,
                                  text badly centered,
                                  inner sep=0pt,
                                 },
                    arrow/.style={draw=black,
                                  very thick,
                                  -latex},
                    bg/.style={rectangle,
                                  draw=black,
                                  rounded corners=1mm,
                                  dashed,
                                  inner sep=1mm,
                                 },
                   ]
\matrix [column sep=2mm, row sep=4mm] {
    \node [     ] (helpS) {};
  & \node [block]  (sys)  {controlled system};
  &
  & \node [block] (input) {calculate new input \eqref{eq:OptControllerGain}, measure $x(k)$, predict $x(k+1)$};\\
  \\\node [     ] (help3) {};\\
    \node [     ] (help2) {};
  & \node [decis] (check) {check for request};
  &
  & \node [decis] (event) {check condition \eqref{eq:triggeringCondition}};
  & \node [     ] (helpB) {};\\
  & \node [block] (solve) {solve QP \eqref{eq:condensedMPC}};
  & \node [     ] (helpX) {~~~~~~};
  & \node [block]  (get)  {request $\gamma$ };\\
  & \node [block] (gamma) {calculate $\gamma$ (see Lem.~\ref{lem:approachesAndEffort})};
  &
  & \node [block]  (law)  {calculate new control law and polytope \eqref{eq:MatricesForAffineLawAndPolytope}};\\
    \node [     ] (help1) {};
  &&&
  & \node [     ] (helpA) {};\\
};
\begin{scope} [every path/.style=arrow]
  \path (check) -- node[near start]{true} (solve);
  \path (solve) -- (gamma);
  \path (gamma) |- (help1.center) -- (help2.center);
  \path (check) -- node[near start,xshift=.6mm]{false} (help2.center) -- (help3.center) -| (check);
  \path (input) -- (event);
  \path (event) -- node[near start]{true} (get) ;
  \path  (get)  --  (law) ;
  \path  (law)  |- (helpA.center) -- (helpB.center);
  \path (event) -- node[near start,xshift=-.5mm]{false} (helpB.center) |- (input);
  \path [rounded corners,dotted]  (get)  -- (helpX.north) |- node[near start,xshift=.5mm,yshift=-4mm]{\begin{rotate}{90}$x(k+1)$ \end{rotate}} (check);
  \path [rounded corners,dotted] (gamma) -| node[near end,xshift=.5mm,yshift=-5mm]{\begin{rotate}{90}$\gamma$ \end{rotate}} (helpX.south) --  (get) ;
  \path [rounded corners,double,transform canvas={yshift=+1.5mm},color=orange] (input) -- node[above]{\color{black}{\; $u^\star(k)$}}  (sys) ;
  \path [rounded corners,double,transform canvas={yshift=-1.5mm},color=orange]  (sys)  -- node[below]{\color{black}{$x(k)$}} (input);
\end{scope}
\begin{pgfonlayer}{background}
  \node [bg, fill=SPS!10, fit=(helpA) (input),xshift=1.6mm,yshift=3.0mm] {}; 
  \node [bg, fill=SPS!25, fit=(helpA) (input),xshift=0.8mm,yshift=1.5mm] {}; 
  \node [bg, fill=HIL!10, fit=(helpS) (sys),xshift=1.6mm,yshift=3.0mm] {}; 
  \node [bg, fill=HIL!25, fit=(helpS) (sys),xshift=0.8mm,yshift=1.5mm] {}; 
  \node [bg, fill=IPC, fit=(help3) (help1) (gamma)] {}; 
  \node [bg, fill=SPS, fit=(helpA) (input)] {}; 
  \node [bg, fill=HIL, fit=(helpS) (sys)] {}; 
\end{pgfonlayer}
\end{tikzpicture}
  \caption{Flow chart of the event-based MPC implementation. Blue box: central node, yellow box: local node, green box: controlled system. Dotted black lines indicate connections that are not used periodically but only in case of an event. Double solid orange lines indicate a permanent feedback control loop.}
	\label{flow:Ablauf}
\end{figure}
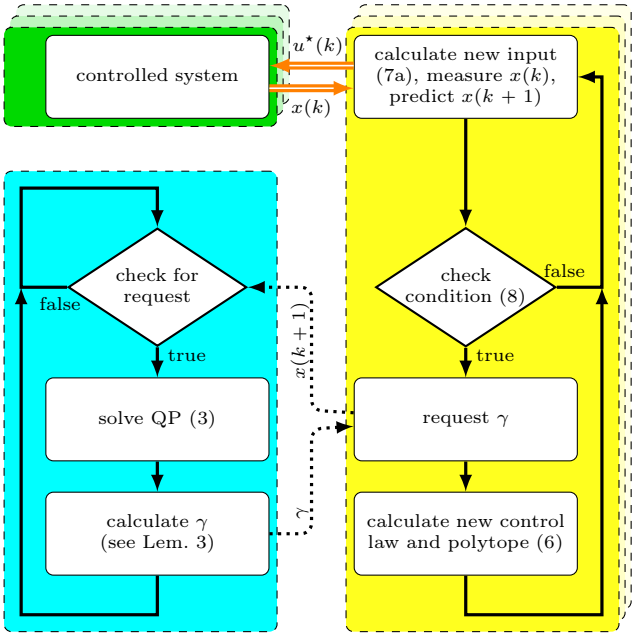%

(ii) If $x(k+1) \not\in \Pe^\star_k $, the QP~\eqref{eq:condensedMPC} must be solved on the central node and a new control law and polytope must be calculated locally.
The time between two consecutive sampling intervals $k$ and $k+1$ can be used to compute the control law and polytope~\eqref{eq:OptControllerGainAndPolytope} that are required in time step $k+1$:
We transfer the predicted state $x(k+1)$ to the central node, solve the quadratic program \eqref{eq:condensedMPC} and transfer data according to approaches~\ref{appr:activeset}-\ref{appr:firstM} back to the local node.
Once the information is available on the local node, matrices of the control law and polytope are calculated.
By solving~\eqref{eq:condensedMPC} for the \textit{predicted state } $x(k+1)$, the control law and polytope are obtained for time step $k+1$.
The control input $u(k+1)$ is calculated at the beginning of the next sampling interval by evaluating the predicted control law. 

The flow chart in Fig.~\ref{flow:Ablauf} visualizes computations and communication in the  HIL setup. In case~(i) all calculations are done on the local node (yellow box). In case~(ii), the central node (blue box) solves the QP~\eqref{eq:condensedMPC}.

\subsection{Validation of the implementation}
We validated the implementation on the PLC and the IPC by regulating the simple double integrator from~\cite{Seron2003} to the origin for 100 random initial states and comparing the state space trajectories that resulted from the HIL setup  to reference trajectories simulated in Matlab. A mean squared error
$\frac{1}{M}\sum\nolimits_{i=1}^M\norm{x(i)_{\rm Matlab}-x(i)_{\rm HIL}} = 0.0071$ resulted, 
where a total of $M= 4051$ time steps arose for the 100 random initial conditions. 
The error results from analog to digital and digital to analog conversion errors, and random noise on the analog signal.

\begin{figure}[b!]
 \centering
 \def\svgwidth{8cm}
 \input{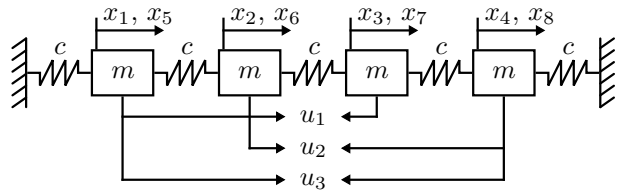}
 \caption[Four-Mass-Oscillator]{Mechanical system with four oscillating masses. $m = c = 1$, no damping. States $x_i,i=1,...,4$ are the position and $x_j,j=5,...,8$ the velocity of the masses. There are three control inputs $u_l,l=1,2,3$.}
 \label{fig:SketchSystemB}
 \end{figure}%

\begin{figure*}[t!]%
\begin{minipage}[t]{0.49 \linewidth}
\centering
%
%
\begin{psfrags}%
\psfragscanon%
%
\psfrag{s05}[t][t]{\color[rgb]{0,0,0}\setlength{\tabcolsep}{0pt}\begin{tabular}{c}$q_\A$\end{tabular}}%
\psfrag{s06}[b][b]{\color[rgb]{0,0,0}\setlength{\tabcolsep}{0pt}\begin{tabular}{c}time in ms \vspace{5px}\end{tabular}}%
\psfrag{s11}[][]{\color[rgb]{0,0,0}\setlength{\tabcolsep}{0pt}\begin{tabular}{c} \end{tabular}}%
\psfrag{s10}[][]{\color[rgb]{0,0,0}\setlength{\tabcolsep}{0pt}\begin{tabular}{c} \end{tabular}}%
\psfrag{s13}[l][l]{\color[rgb]{0,0,0}task a}%
\psfrag{s14}[l][l]{\color[rgb]{0,0,0}\setlength{\tabcolsep}{0pt}\begin{tabular}{l}\vspace{-2px} task b\,($\Delta t_{\rm b}$)\end{tabular}}%
\psfrag{s15}[l][l]{\color[rgb]{0,0,0}task c}%

%
\psfrag{x01}[t][t]{$0$}%
\psfrag{x02}[t][t]{$5$}%
\psfrag{x03}[t][t]{$10$}%
\psfrag{x04}[t][t]{$15$}%
\psfrag{x05}[t][t]{$20$}%
\psfrag{x06}[t][t]{$25$}%
\psfrag{x07}[t][t]{$30$}%
%
\psfrag{v01}[r][r]{$0$}%
\psfrag{v02}[r][r]{$100$}%
\psfrag{v03}[r][r]{$200$}%
\psfrag{v04}[r][r]{$300$}%
\psfrag{v05}[r][r]{$400$}%
\psfrag{v06}[r][r]{$500$}%
\psfrag{v07}[r][r]{$600$}%
%
\begin{tikzpicture}
    \node[anchor=south west,inner sep=0] at (0,0) {\includegraphics[width=8cm, height=7cm]{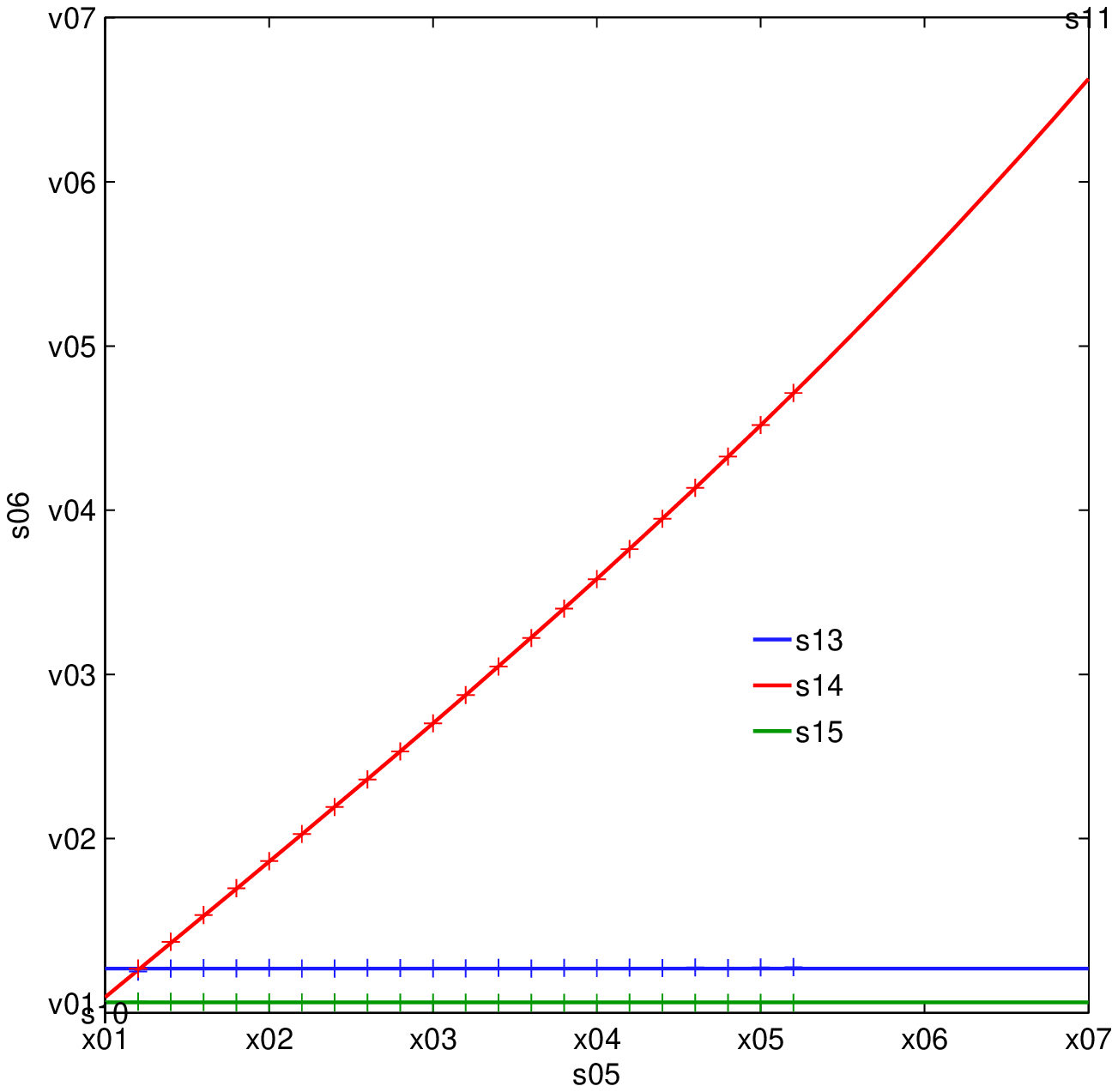}};
    \draw[dashed,red,thick] (0.79,0.55) rectangle (1.3,1.1);
    \begin{scope}[shift={(0,0)}]
    %
    \scriptsize
\psfrag{x01}[t][t]{$0$}%
\psfrag{x02}[t][t]{$1$}%
\psfrag{x03}[t][t]{$2$}%
%
\psfrag{v01}[r][r]{$0$}%
\psfrag{v02}[r][r]{}%
\psfrag{v03}[r][r]{$20$}%
\psfrag{v04}[r][r]{}%
\psfrag{v05}[r][r]{$40$}%
\psfrag{v06}[r][r]{}%
    	%
    	
    	\draw[dashed,red,thick,fill=white] (1.5,4.5) rectangle (4.5,6.6);
    	\draw[red, thick, dotted] (0.79,1.1) -- (1.5,4.5);
    	\draw[red, thick, dotted] (1.3,1.1) -- (4.5,4.5);
    	\node[anchor=south west,inner sep=0] at (1.8,4.65) {\includegraphics[width = 2.45cm]{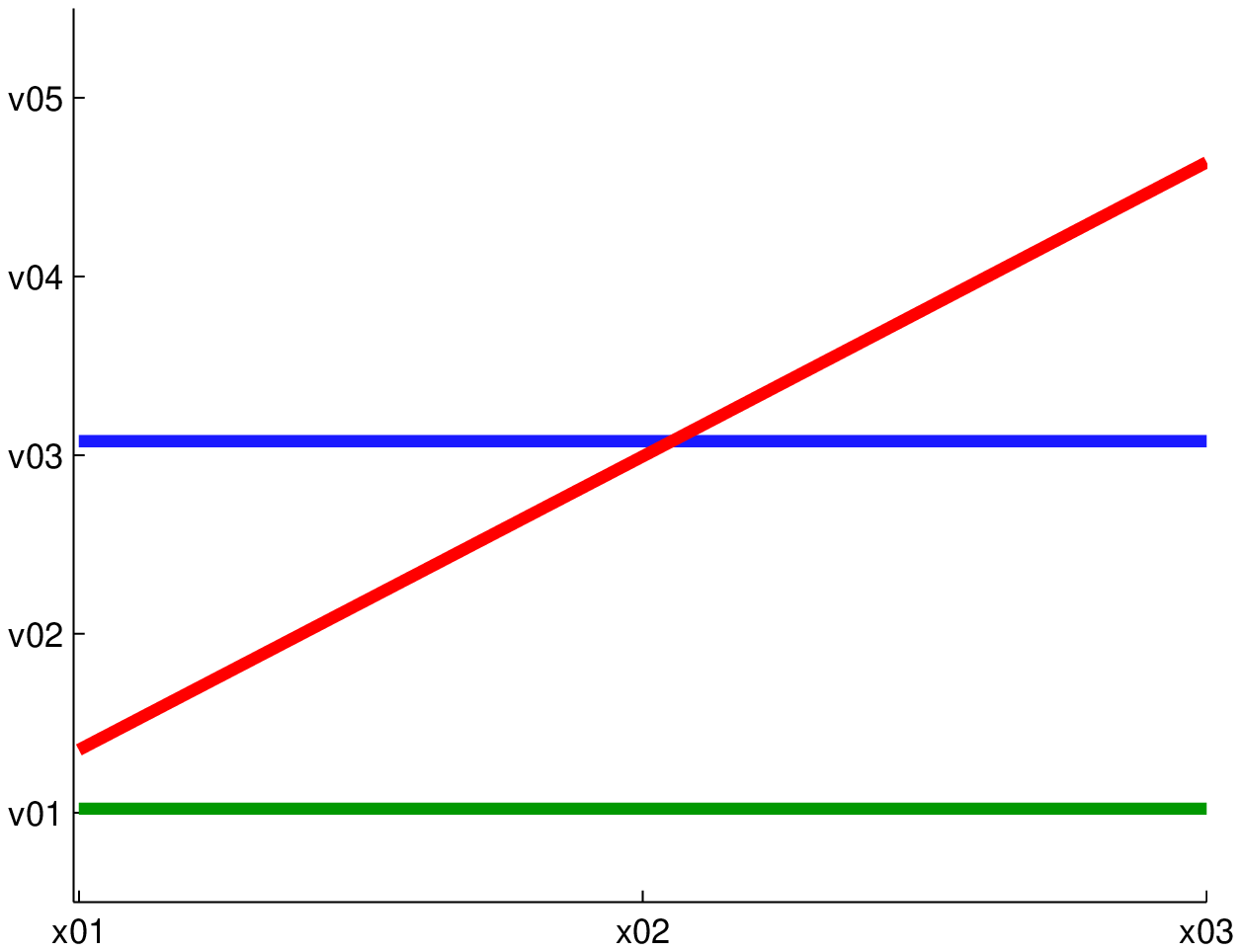}};
    	    	
    \end{scope}
\end{tikzpicture}
\end{psfrags}%
%

  \caption{Comparison of computation times for tasks~a-c on the local node (PLC) as a function of $q_\A$, which is bounded from above by $30$  
           according to~\eqref{eq:UpperBoundActiveConstraints}. The random initial conditions used in the experimental study cover the cases $q_\A\in\{0, \dots, 21\}$. All measurements were repeated ten times. The standard deviations, which amount to less than $2.6\%$, $0.2\%$ and $1.5\%$ for tasks a, b and c, respectively, are not shown, because most of the resulting error bars would hardly be visible.}
  \label{fig:Times2}
\end{minipage}
\hfill
\begin{minipage}[t]{0.49 \linewidth}
\centering
%
%
\begin{psfrags}%
\psfragscanon%
%
\psfrag{s05}[l][l]{\color[rgb]{0,0,0}\setlength{\tabcolsep}{0pt}\begin{tabular}{l}  \,$\Delta t_{\rm mat}$\end{tabular}}%
\psfrag{s06}[l][l]{\color[rgb]{0,0,0}\setlength{\tabcolsep}{0pt}\begin{tabular}{l}  \,$\Delta t_{\rm inv}$\end{tabular}}%
\psfrag{s08}[t][t]{\color[rgb]{0,0,0}\setlength{\tabcolsep}{0pt}\begin{tabular}{c}$q_\A$\end{tabular}}%
\psfrag{s09}[b][b]{\color[rgb]{0,0,0}\setlength{\tabcolsep}{0pt}\begin{tabular}{c}time in ms \vspace{5px}\end{tabular}}%
\psfrag{s12}[][]{\color[rgb]{0,0,0}\setlength{\tabcolsep}{0pt}\begin{tabular}{c} \end{tabular}}%
\psfrag{s13}[][]{\color[rgb]{0,0,0}\setlength{\tabcolsep}{0pt}\begin{tabular}{c} \end{tabular}}%
\psfrag{s14}[l][l]{}%
\psfrag{s15}[l][l]{\color[rgb]{0,0,0}task b ($\Delta t_{\rm b} = \Delta t_{\rm inv} + \Delta t_{\rm mat} $)}%
%
\psfrag{x12}[t][t]{$0$}%
\psfrag{x13}[t][t]{$5$}%
\psfrag{x14}[t][t]{$10$}%
\psfrag{x15}[t][t]{$15$}%
\psfrag{x16}[t][t]{$20$}%
\psfrag{x17}[t][t]{$25$}%
\psfrag{x18}[t][t]{$30$}%
%
\psfrag{v12}[r][r]{$0$}%
\psfrag{v13}[r][r]{$100$}%
\psfrag{v14}[r][r]{$200$}%
\psfrag{v15}[r][r]{$300$}%
\psfrag{v16}[r][r]{$400$}%
\psfrag{v17}[r][r]{$500$}%
\psfrag{v18}[r][r]{$600$}%
%
\includegraphics[width=8cm, height=7cm]{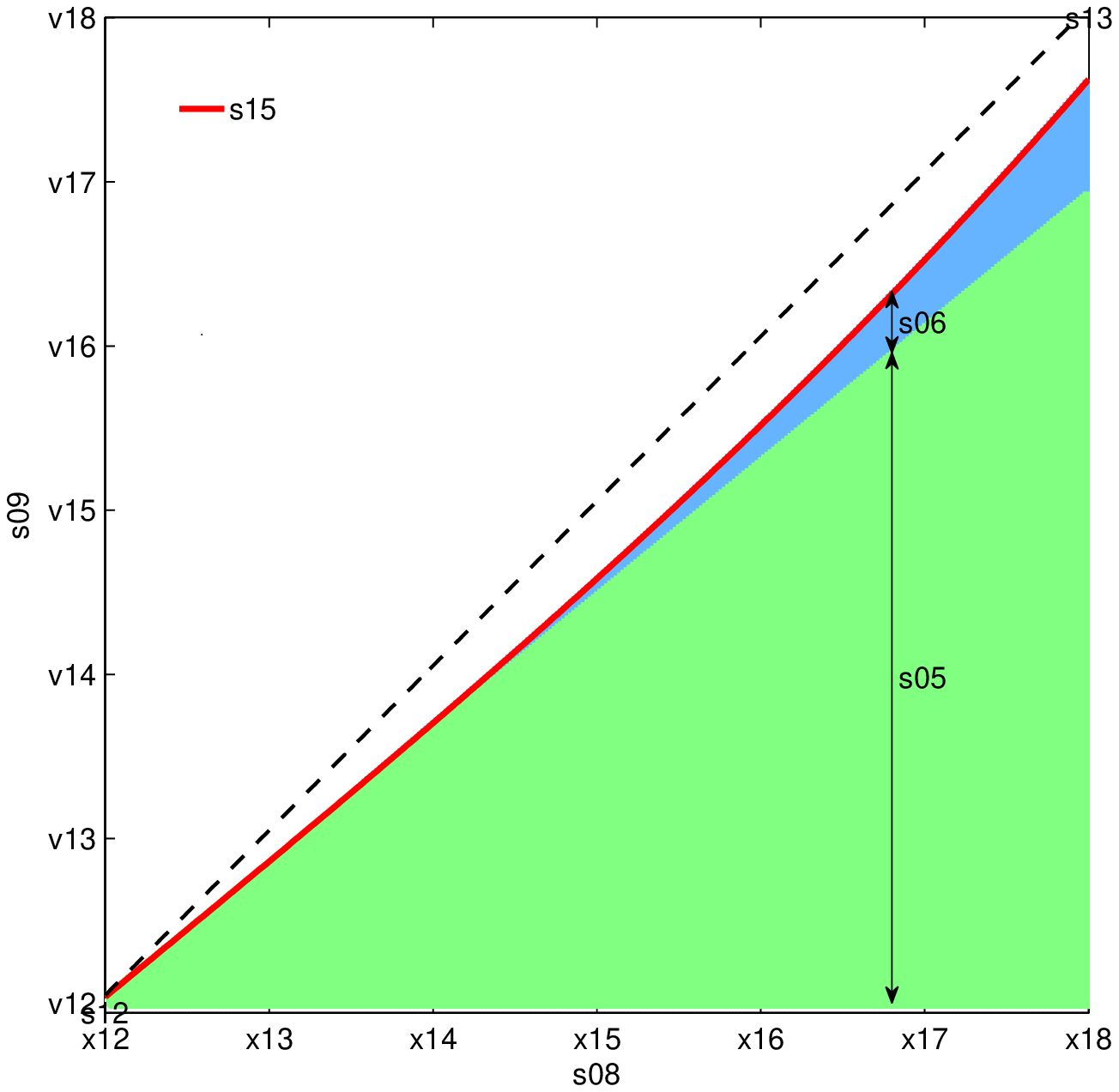}%
\end{psfrags}%
%

  \caption{Contributions to $\Delta t_b= \Delta t_{{\rm inv}}+ \Delta t_{{\rm mat}}$. The dashed line marks the upper bound stated in~\eqref{eq:UpperBoundtb}. 
}
  \label{fig:Times}
\end{minipage}
\end{figure*}%

\section{Results}\label{sec:results}
We measure computation times of approach~\ref{appr:activeset} running on the HIL setup described in Section~\ref{sec:implementation}. The results corroborate Proposition~\ref{prop:upperBoundKappa}.

\subsection{Example system}
\label{sec:example}
Consider the mechanical system with four oscillating masses connected by springs sketched in Fig.~\ref{fig:SketchSystemB}. 
Three inputs ($u_l,l=1,2,3$) are used to control tension between the masses. Masses and spring constants have the value one and there is no damping. The resulting continuous time  system
%
\begin{align}
\dot{x}&= Ax+Bu \text{, with} \label{eq:system}\\
\scalebox{0.99}{$A$} &\scalebox{0.9}{$= \begin{bmatrix}
  0^{4\times4} & I^{4\times4}\\
  -F_c          & 0^{4\times4}
\end{bmatrix} \text{,  }
B = \begin{bmatrix}
0^{4\times3} \\
F_u
\end{bmatrix} \text{ ,}$} \nonumber \\
\scalebox{0.99}{$F_c$} &\scalebox{0.9}{$= \begin{bmatrix}
2 &-1 & 0 & 0\\
-1 & 2 &-1 & 0\\
0 &-1 & 2 &-1\\
0 & 0 &-1 & 2
\end{bmatrix}
~\text{and}~F_u = \begin{bmatrix}
1 & 0 & 1 \\
0 & 1 & 0 \\
-1 & 0 & 0 \\
0 &-1 &-1
\end{bmatrix}$} \nonumber
\end{align}
%
is discretized with zero order hold and a sampling time of $T_s=0.5$ seconds. State variables $x\in\R^8$ and input variables $u\in\R^3$ are subject to the constraints $x_i \in \left[-4,4\right], i=1,...,8$ and $u_l\in \left[-0.5,0.5\right], l=1,2,3$.
Cost function matrices are chosen to be $Q=I^{n\times n}$ and $R=I^{m\times m}$, and $P$ is set to the solution of the discrete-time algebraic Riccati equation. The horizon of the MPC problem is chosen to be $N=10$. 
The number of simultaneously active constraints is bounded above according to
\begin{equation}\label{eq:UpperBoundActiveConstraints}
  q_\A \leq mN = 30.
\end{equation}

\subsection{Hardware in the loop test}\label{sec:HILresults}
Approach~\ref{appr:activeset} provides the best trade-off of  data transmission to computational cost on the local node. 
We conduct measurements of the computation times on the local node with the HIL setup described in Section~\ref{sec:implementation} to analyze approach~\ref{appr:activeset} in more detail. 
Specifically, we report actual computational times of the HIL system for the following tasks to corroborate Proposition~\ref{prop:upperBoundKappa}:
\begin{itemize}
\item[a)] transferring the current system state to the central node, solving the QP on the central node, and sending $\gamma$ as defined in~\eqref{eq:bitwiseRepresentationGamma} to the local node,
\item[b)] evaluation of~\eqref{eq:MatricesForAffineLawAndPolytope} on the local node,
\item[c)] computing the control input and evaluating  condition~\eqref{eq:OptControllerGainAndPolytope} on the local node.
\end{itemize}
We assume that flops and computation times are related by a constant factor, i.e. 
\begin{equation}\label{eq:DeltaPropEta}
  \Delta t_{\rm inv} = \alpha \eta_{\rm inv}, \quad
  \Delta t_{\rm mat} = \alpha \eta_{\rm mat}
\end{equation}
for some $\alpha \in \R^+$. 

The computational times that result for 100 random feasible initial conditions are summarized in Fig.~\ref{fig:Times2}. 
The times required for tasks~a, b and c can be approximated by a constant, a cubic function, and a constant, respectively. 
The time required for task~c is obviously small compared to those of tasks~a and~b, which demonstrates that the control input can be computed with very low computational effort if $e(x)=0$. The time required for task~b dominates the computational time required on the local node for all but very small $q_\A$.

As suggested by Proposition~\ref{prop:upperBoundKappa}, we further analyze the computation time $\Delta t_b$ for task~b by splitting it into the time $\Delta t_{{\rm inv}}$ necessary to calculate the inverse of $G_{\A}H^{-1}(G_\A)'$ and the time $\Delta t_{{\rm mat}}$ required for all remaining matrix operations in~\eqref{eq:MatricesForAffineLawAndPolytope}. 
Figure~\ref{fig:Times} shows the two contributions 
\begin{align}\label{eq:sumTinvTmat}
\Delta t_b = \Delta t_{{\rm inv}}+\Delta t_{{\rm mat}} .
\end{align}
Substituting  
\begin{equation}\label{eq:ActualKappa}
  \frac{\Delta t_{{\rm inv}}}{\Delta t_{{\rm mat}}}\le \frac{18}{79}\approx 0.23\,,
\end{equation}
which holds according to Proposition~\ref{prop:upperBoundKappa} and with~\eqref{eq:DeltaPropEta},
results in 
\begin{equation}\label{eq:UpperBoundtb}
  \Delta t_{b} \le \left( 1+ \frac{18}{79} \right) \Delta t_{{\rm mat}} \,.
\end{equation}
This can be understood as an upper bound on $\Delta t_{b}$, which is also shown in Fig.~\ref{fig:Times}.
It is evident from Fig.~\ref{fig:Times} that the absolute contribution of $\Delta t_{{\rm inv}}$ remains small despite the cubic dependence of $\Delta t_{{\rm inv}}$ on $q_\A^3$. 
We conclude approach~\ref{appr:activeset} is the best one of the four variants, because the inversion of $G_{\A}H^{-1}(G_\A)'$ required in~\ref{appr:activeset} but not in~\ref{appr:sendInverse} results in only a small additional computational effort on the local node. Note that this result does not only apply to the example, but the bound~\eqref{eq:kappa} established in Proposition~\ref{prop:upperBoundKappa} applies in general. 

\section{Conclusion and outlook}\label{sec:conclusion}
We analyzed the practical usefulness of a networked event-triggered MPC approach by implementing it on industrial hardware. 
The proposed approach combines a powerful compute node, which solves optimal control problems on demand, with lean local nodes that control local systems with simple affine, regional optimal feedback laws. The local nodes use their feedback laws as long as they are optimal and otherwise request a new law from the central node. Several options for an implementation exist that differ with respect to their balance of local computational effort and network bandwidth requirements. We extended earlier theoretical results to determine the best variant among these choices and corroborated the theoretically expected features of the selected variant with an implementation on industrial control hardware.

An extension to tube-based robust MPC is possible (see~\cite{SchulzeDarup2017} for a simulation-based analysis). We believe an extension in this direction is more of theoretical interest, because the case treated in the present paper inherits the robustness of the regional approach~\cite{BernerP2016}.  This aspect is currently under investigation.

The proposed approach is based on the idea that once the set of active constraints is known for the current system state, an optimal feedback law and the polytope on which it is optimal are known (see Lemma~\ref{lem:OptionsForTransmittedData}). This also holds in nonlinear MPC with discrete-time system models~\cite{Monnigmann2015,Dominguez2011}. 
It therefore suffices to send the active sets across the network in nonlinear MPC. However, the local nodes must solve nonlinear systems of equations in this case, while affine control laws only need to be evaluated in the linear MPC case. It remains to investigate how the computational effort on the local nodes can be reduced, for example by sending redundant information as in approach~\ref{appr:sendInverse}.  

Finally, future work must combine the proposed approach with methods for predicting the sequence of affine laws along the closed-loop trajectory~\cite{BernerKoenig2019}.
This will help reducing the number of optimal control problems further in particular far from the origin, where the affine laws are not reused frequently.

\section*{Acknowledgment}
Support by the \textit{Deutsche Forschungsgemeinschaft} (DFG) under the grant MO~1086/15-1 is gratefully acknowledged.

\bibliographystyle{abbrv}
\bibliography{ComplexityAnalysis}

\end{document}